\documentclass[12pt,reqno]{amsart}

\usepackage[stable]{footmisc}

\numberwithin{equation}{section}

\usepackage{mdwlist}

\usepackage{nameref}

\usepackage[final, 
color]{showkeys}
\definecolor{refkey}{gray}{.85}
\definecolor{labelkey}{gray}{.85}

\makeatletter
\let\orgdescriptionlabel\descriptionlabel
\renewcommand*{\descriptionlabel}[1]{%
  \let\orglabel\label
  \let\label\@gobble
  \phantomsection
  \edef\@currentlabel{#1}%
  \let\label\orglabel
  \orgdescriptionlabel{#1}%
}
\makeatother

\usepackage{undertilde}

\usepackage{AKstyle}

\usepackage{caption}
\usepackage[labelfont=rm]{subcaption}

\usepackage[pagebackref=true, colorlinks]{hyperref}

\hypersetup{pdffitwindow=true,linkcolor=blue,citecolor=blue,urlcolor=blue,menucolor=blue}

\usepackage{comment}

\begin{document}

\author{Alex Kontorovich}
\thanks{
The author gratefully acknowledges support from
an NSF CAREER grant DMS-1254788, an Alfred P. Sloan Research Fellowship, a Yale Junior Faculty Fellowship, and support at IAS from The Fund for Math and The Simonyi Fund.}
\email{alex.kontorovich@yale.edu}
\address{Math Department, Yale University, New Haven, CT 06511 USA, and School of Mathematics, Institute for Advanced Study, Princeton, NJ 08540 USA}

\title
{Levels of Distribution and the Affine Sieve}

\date{\today}
\maketitle
\setcounter{tocdepth}{2}
\tableofcontents

\section{Introduction}\label{sec:intro}

This article is an expanded version of the author's lecture in the
Basic Notions Seminar
at Harvard
, September 2013. 
Our 
goal is a brief 
and
introductory
exposition of 
aspects of two topics in sieve theory which have received 
attention recently: (1) the spectacular work of Yitang 
Zhang, under the title ``Level of Distribution,'' and (2) 
the so-called ``Affine Sieve,''
introduced by Bourgain-Gamburd-Sarnak.

\section{Level of Distribution for the Primes}

Let $p_{n}$ be the $n$th prime number. 
We begin with the infamous 

\vskip.1in
\noindent
{\bf Twin Prime Conjecture:} 
$$
\liminf_{n\to\infty}( p_{n+1}-p_{n})\ =\ 2.
$$

A slight weakening of this problem is called the 

\phantomsection
\vskip.1in
\noindent
{\bf Bounded Gaps Conjecture:}\label{conj:BddGaps} 
$$
\liminf_{n\to\infty} (p_{n+1}-p_{n})\ <\ \infty.
$$

A tremendous shock ran through the mathematical community in April 2013 when 
Yitang
Zhang \cite{Zhang2013} proved 

\noindent
\phantomsection
{\bf Zhang's Theorem (2013):}\label{thm:Zhang} 
{\it The Bounded Gaps Conjecture is true. In particular,} 
$$
\liminf_{n\to\infty} (p_{n+1}-p_{n})\ <\ 7\times10^{7}.
$$


Our goal in this section is to explain what is meant by a ``level of distribution'' for the primes, and give some hints of the role it plays in the proof of Zhang's theorem.

\subsection{The Distribution of Primes}\

First we recall the Prime Number Theorem (PNT), 
proved independently and simultaneously by Hadamard
 \cite{Hadamard1896}
 and 
 de la Vall\'ee Poussin
  \cite{delaValleePoussin1896}
 in 1896, following the strategy introduced in Riemann's 1859 epoch-making memoir  \cite{Riemann1859}. 
It is often stated as: 
$$
\pi(x)\ :=\ \sum_{p<x}1\quad\sim\quad {x\over \log x},\qquad\qquad 
x\to\infty,
$$
where, as throughout, $\log$ is to base $e$, and $p$ denotes a prime.
The first Basic Notion is that this is the ``wrong'' formula, not in the sense of being untrue, but in the sense that
\be\label{eq:PNTbad}
\pi(x)\ =\ {x\over \log x} +\gW\left({x\over \log ^{2}x}\right)
,
\ee
the error term being unnecessarily large. (Here $\gW$ is the negation of little-oh.) A more precise statement 
of 
PNT, not far from the best currently known,  is the following.

\phantomsection
\vskip.1in
\noindent
{\bf Prime Number Theorem:}\label{thm:PNT}
{\it For any $A>1$, }
\be\label{eq:PNT}
\pi(x)\ =\ \Li(x)+O_{A}\left({x\over \log ^{A}x}\right)
,\qquad\qquad  \text{ as } x\to\infty.
\ee

Here the subscript $A$ in the big-Oh means that the implied constant 
depends on $A$, and $\Li$ is the ``logarithmic integral'' function
$$
\Li(x)\ :=\ \int_{2}^{x}{dt\over\log t}
.
$$
By an exercise in partial integration, we have that
$$
\Li(x)
\ =\
{x\over \log x}+{x\over \log^{2}x} +O\left({x\over \log^{3}x}\right)
,
$$
which 
together with \eqref{eq:PNT} 
implies \eqref{eq:PNTbad}.
On the other hand, the Riemann Hypothesis (RH) predicts that
$$
\pi(x)\ =\ \Li(x)+O(\sqrt x\log x),
$$
epitomizing the ``square-root cancellation'' phenomenon. If true, this estimate would be  best possible (up to $\log$ factors), as Littlewood proved in 1914 that
$$
\pi(x)\ =\ \Li(x)+\gW\left(\sqrt x\ { \log\log\log x\over\log x}\right).
$$
In fact, he showed that the difference $\pi(x)-\Li(x)$ infinitely-often attains both positive and negative values of this order of magnitude.

\subsection{Primes in Progressions}\

The next most basic question is: How are the primes distributed in arithmetic progressions? 
Given an integer $q\ge1$, often called the ``level'' 
in this context, and a coprime number $(a,q)=1$, let
$$
\pi(x;a,q)\ :=\ \sum_{p<x\atop p\equiv a(\mod q)}1
$$
denote the number of primes up to $x$ in the progression $a(\mod q)$. A relatively minor modification to the proof of \eqref{eq:PNT} gives

\vskip.1in
\noindent
{\bf PNT in Progressions:} {\it For any $A>1$,}
\be\label{eq:PNTAP}
\pi(x;a,q)\ =\ {\Li(x)\over \phi(q)}+O_{A,q}\left({x\over \log^{A}x}\right),\qquad\qquad x\to\infty. 
\ee

Meanwhile, the 
Generalized
Riemann Hypothesis (GRH) predicts
\be\label{eq:APsRH}
\pi(x;a,q)\ =\ {\Li(x)\over \phi(q)}+O_{\vep}\left(x^{1/2+\vep}\right),
\ee
for any $\vep>0$.
These estimates
confirm our intuition that primes should not favor one primitive  (
meaning $a$ and $q$ are coprime%
) arithmetic progression $\mod q$ over others, there being $\phi(q)=|(\Z/q\Z)^{\times}|$ of them total.

In applications, it is often important to be able to use formulae like \eqref{eq:PNTAP} while allowing $q$ to vary with $x$.
For example, 
it does not seem unreasonable that we should be able to  use \eqref{eq:PNTAP} to estimate, 
say, the number of primes up to $e^{100}$ which 
are 
$1$ mod $
100^{3}$, or primes up to $e^{1000}$ which are $1$ mod $1000^{3}$, or more generally, primes up to $x=e^{\ell}$ which are $1\mod \ell^{3}$. 
In these examples, we have taken a level of size
 $q=\ell^{3}=\log^{3}x$ 
 which, it turns out, is growing too rapidly relative to $x$ to obtain a meaningful asymptotic from present methods;
the error terms in all of these questions might swamp the main terms, giving no estimate at all. 

To address this issue of uniformity in the level $q$, there is a famous 
estimate proved by  Walfisz \cite{Walfisz1936}  in 1936 by adapting work of Siegel   \cite{Siegel1935a}.

\phantomsection
\vskip.075in
\noindent
{\bf Siegel-Walfisz Theorem:}\label{thm:SiegWalf}
{\it Given any positive constants $A$ and $B$, any $q<\log ^{B}x$, and any $(a,q)=1$, we have}
\be\label{eq:SW}
\pi(x;a,q)\ =\ {\Li(x)\over \phi(q)}+O_{A,B}\left({x\over \log ^{A}x}\right).
\ee

It may appear that the uniformity issue 
in the range $q<\log^{B}x$ has been completely resolved, but there's a catch: 
the implied constant in 
\eqref{eq:SW}
 coming from the proof is  ``ineffective.'' This means that, once  the 
 parameters
 $A$ and $B$ are supplied, there is no known procedure to 
 determine
  the constant. 
  Thus we still have no way of
verifying that
$\Li(e^{\ell})/\phi(\ell^{3})$ is
an accurate estimate
 for $\pi(e^{\ell};1,\ell^{3})$.
This so-called ``Siegel zero'' phenomenon is the sense in which we do not know the PNT in progressions.

The danger of an ineffective constant is beautifully 
illustrated
by 
 Iwaniec's (facetious)

\vskip.075in
\noindent
{\bf Theorem:}
{\it There exists a constant $C>0$ such that, if RH holds up to height $C$ (meaning $\gz(\gs+it)\neq0$ for all $\foh<\gs<1$, $|t|<C$), then RH holds everywhere.}

\vskip.1in
This fantastic result seems to reduce RH to a finite computation; before we get too excited, let's have a look at the
\pf
There are two cases. 

{\bf Case 1:} Assume RH is true. Set $C=1$, and RH holds.

{\bf Case 2:} Assume RH is false, that is, $\gz(\gs+it)=0$ for some $\foh<\gs<1$ and some $t>0$. Set $C=t+1$. 
The statement is vacuously true.
\epf

As an aside, we briefly recall that a similar phenomenon occurs in the study of Gauss's Class Number Problem. Let $-d<0$ be 
the discriminant of an imaginary quadratic field and let  $h(-d)$ be the corresponding class number
 (see wikipedia for definitions, which will not be needed for our discussion). In 1936, Siegel (based on earlier work by Hecke \cite{Landau1918},   Deuring \cite{Deuring1933}, Mordell \cite{Mordell1934}, Heilbronn \cite{Heilbronn1934}, and Landau \cite{Landau1935}) proved that
\be\label{eq:ClassN}
h(-d)\ \gg_{\vep}\ d^{1/2-\vep},
\ee
for any $\vep>0$. Again this implied constant is ineffective, and thus does not allow one to, e.g., tabulate all $d$ with $h(-d)=1$ (the Class Number One Problem). 
Much later, Goldfeld \cite{Goldfeld1976, Goldfeld1985}
 (1976) together with Gross-Zagier \cite{GrossZagier1986} (1985) managed to circumvent the ineffectivity, proving
\be\label{eq:ClassEff}
h(-d)\ >\ \frac1{55}\log d,
\ee
whenever $d$ is prime 
 (we make this restriction only to give the simplest formula). Thanks to \eqref{eq:ClassEff} (and much other work), we now have complete tables of all $d$ with class number up to $100$.
The point of this aside is that, just because one proof gives an ineffective constant, there might be a completely different proof for which the constants are absolute. 
Resolving this 
``Siegel zero'' issue is one of the main outstanding problems in analytic number theory. 

\subsection{Primes in Progressions on Average}\

In many applications, what is needed is not uniformity for a single level $q$, but over a range of $q$. 
This is the heart of what is meant by a ``level of distribution,'' as explained below.

Assuming GRH, we see from \eqref{eq:APsRH} that 
\be\label{eq:GRHlevel}
\sum_{q<Q}\max_{(a,q)=1}\left|\pi(x;a,q)-{\Li(x)\over \phi(q)}\right|
\ \ll_{\vep}\
\sum_{q<Q}x^{1/2+\vep}
\ <\
Qx^{1/2+\vep}
.
\ee
So if we take $Q=x^{1/2-2\vep}
$, say, then the error terms add up to at most $
x^{1-\vep}$, while there are about $x/\log x$ primes up to $x$. That is, all of these errors summed together still do not exceed the total number of primes.
This immediately leads us to the

\noindent
{\bf Definition: Level of Distribution} (for primes in progressions). We will say that the primes have a {\it level of distribution} $Q$ if,
for all $A<\infty$,
%
\be\label{eq:LevelDef}
\sum_{q<Q}\max_{(a,q)=1}\left|\pi(x;a,q)-{\Li(x)\over \phi(q)}\right|
\ =\
O_{A}\left({x\over \log^{A} x}\right)
.
\ee
When $Q$ can be taken as large as $x^{\vt-\vep}$ for some $\vt>0$, 
we call $\vt$ an {\it exponent of distribution} for the primes.

Note that 
level of distribution is not a quantity inherent to the sequence of primes, but is instead a function of what one can {\it prove} about the primes.
While GRH implies the level $Q=x^{1/2-\vep}$ (or exponent $\vt=1/2$),  the unconditional Siegel-Walfisz estimate \eqref{eq:SW} gives only a level of size $Q=\log^{A}x$, which is not even a positive exponent $\vt$. 

It was a dramatic breakthrough when  Bombieri \cite{Bombieri1965}
 and A. I. Vinogradov \cite{Vinogradov1965} (based on earlier work of Linnik \cite{Linnik1941}, Renyi \cite{Renyi1948}, Roth \cite{Roth1965}, and Barban \cite{Barban1966}) independently and simultaneously proved the

\phantomsection
\vskip.1in
\noindent
{\bf Bombieri-Vinogradov Theorem (1965):}\label{thm:BV} {\it The primes have 
exponent of distribution $\vt=1/2$. More precisely, for any constant $A>1$, there exists a constant $B>1$ so that}
$$
\sum_{q<{x^{1/2}\over \log^{B} x}}\max_{(a,q)=1}\left|\pi(x;a,q)-{\Li(x)\over \phi(q)}\right|
\ \ll_{A}\
{x\over \log ^{A}x}
.
$$

The Bombieri-Vinogradov theorem (B-V) is thus an unconditional substitute for GRH on average, 
since both produce
the same exponent of distribution $\vt=1/2$!
(The implied constant is still ineffective, as the proof uses \hyperref[thm:SeigWalf]{Siegel-Walfisz}; we have not escaped the ``Siegel zero'' problem.)

Being even more ambitious
, one may ask for variation in the size of the error term; after all, we crudely imported the worst possible error from \eqref{eq:APsRH}
into \eqref{eq:GRHlevel}. Applying a ``square-root cancellation'' philosophy yet again, one might boldly posit that the term $Qx^{1/2}$ on the right side of \eqref{eq:GRHlevel} can be replaced by $Q^{1/2}x^{1/2}$, in which case 
$Q$ can be taken as large as $x^{1-\vep}$.
This is the

\phantomsection
\vskip.1in
\noindent
{\bf Elliott-Halberstam Conjecture  \cite{ElliottHalberstam1968} (1968):}\label{conj:EH} {\it The primes have exponent of distribution $\vt=1$. That is, for any $\vep>0$ and $A<\infty$,}
\be\label{eq:EH}
\sum_{q<x^{1-\vep}}\max_{(a,q)=1}\left|\pi(x;a,q)-{\Li(x)\over \phi(q)}\right|
\ =\
O_{A,\vep}\left({x\over \log^{A} x}\right)
.
\ee

The Elliott-Halberstam Conjecture (E-H), if true, goes far beyond any RH-type statement, as far as we are aware. 
As long as we are already dreaming, we may as well suppose that this further square-root cancellation happens not only on average, as E-H claims, but individually; by this we mean the following. Returning to \eqref{eq:APsRH},
the ``main'' term is very roughly of size $x/q$, so might not the error be of square-root the main term, not just square-root of $x$? This is 

\vskip.1in
\noindent
{\bf Montgomery's Conjecture
\cite{Montgomery1971}
 (1971):} For all $\vep>0$,
$$
\pi(x;a,q)\ =\ {\Li(x)\over \phi(q)}+O_{\vep}\left({x^{1/2+\vep}\over q^{1/2}}\right)
.
$$

Montgomery's Conjecture immediately implies E-H, but
we emphasize again that both of these assertions are not, as far as we know, consequences of any RH-type statement. 

Nothing beyond \hyperref[thm:BV]{B-V} has ever been proved towards the pure level of distribution defined in \eqref{eq:LevelDef}. But if one drops the absolute values, fixes one non-zero integer $a$, and weights the errors at level $q$ by a function $\gl(q)$ which is ``well-factorable'' (the precise meaning of which we shall not give here), then one can go a bit into the \hyperref[conj:EH]{E-H} range. Building on work by Fouvry-Iwaniec \cite{FouvryIwaniec1983, Fouvry1984}, we have the

\phantomsection
\vskip.1in
\noindent
{\bf Bombieri-Friedlander-Iwaniec Theorem \cite{BombieriFriedlanderIwaniec1986} (1986):}\label{thm:BFI} {\it Fix any $a\neq0$ and let $\gl(q)$ be a  ``well-factorable'' function. Then for any $A>1$ and $\vep>0$,}
\be\label{eq:BFI}
\sum_{q<x^{4/7-\vep}}\gl(q)\left(\pi(x;a,q)-{\Li(x)\over \phi(q)}\right)
\ \ll_{a,A,\vep}\ 
{x\over \log ^{A}x}
.
\ee

Thus in the weighted sense above, the Bombieri-Friedlander-Iwaniec Theorem (BFI) gives a weighted exponent of distribution 
$$
\vt\ = \ 4/7 \quad >\quad 1/2
,
$$ 
giving some partial evidence towards \hyperref[conj:EH]{E-H}. But before we get too optimistic about the full E-H, let us point out just how delicate the conjecture is. Building on 
work of Maier
\cite{Maier1985}, Friedlander and Granville
\cite{FriedlanderGranville1989}
 showed that the level $x^{1-\vep}$ in \eqref{eq:EH} cannot be replaced by $x(\log x)^{-A}$. More precisely, we have the following

\vskip.1in
\noindent
{\bf Friedlander-Granville Theorem (1989):} {\it For any $A>0$, there exist arbitrarily large values of $a$ and $x$ for which}
$$
\sum_{q<x(\log x)^{-A}\atop (q,a)=1}\left|\pi(x;a,q)-{\Li(x)\over \phi(q)}\right|
\ \gg_{A}\ 
{x\over \log x}
.
$$

In particular, the asymptotic formula $\pi(x;a,q)\sim \Li(x)/\phi(q)$ can be false for $q$ as large as $x/\log^{A}x$.

\subsection{Small Gaps Between Primes}\

Let us return to the \hyperref[conj:BddGaps]{Bounded Gaps Conjecture} (now Zhang's Theorem).
Recall that $p_{n}$ is the $n$th prime. Before studying absolute gaps, what can we say about gaps relative to the average? \hyperref[thm:PNT]{PNT} tells us that $p_{n}\sim n\log n$, so the average gap $p_{n+1}-p_{n}$ is of size about $
\log p_{n}$. Hence 
\be\label{eq:AvgGap}
\gD\ := \ \liminf_{n\to\infty} {p_{n+1}-p_{n}\over \log p_{n}}\ \le\  1.
\ee
Here is a very abbreviated history on reducing the number on the right side of \eqref{eq:AvgGap}.
\\

\hskip-.5in
\begin{tabular}{|l|ccc|l|}
\hline
Hardy-Littlewood \cite{Rankin1940} (1926): &$\gD$&$\le$&$2/3$, &\ assuming GRH.\\
\hline
Erd\"os \cite{Erdos1940} (1940): &$\gD$&$<$&$1$, &\ unconditional; by sieving.\\
\hline
Bombieri-Davenport \cite{BombieriDavenport1966} (1966): &$\gD$&$\le$&$1/2$, &\begin{tabular}{l}unconditional; by refining\\ Hardy-Littlewood and\\replacing GRH by \hyperref[thm:BV]{B-V}.\end{tabular}\\
\hline
Maier  \cite{Huxley1977, Maier1988} (1988): &$\gD$&$<$&$1/4$, &\begin{tabular}{l}unconditional; using a \\radically different method 
\end{tabular}\\
\hline
Goldston \cite{Goldston1992} (1992): &$\gD$&$=$&$0$, &\begin{tabular}{l}assuming \hyperref[conj:EH]{E-H} and another\\ E-H type conjecture.\end{tabular}\\
\hline
Goldston-Pintz-Y{\i}ld{\i}r{\i}m \cite{GoldstonPintzYildirim2009} (2005): &$\gD$&$=$&$0$, &\ unconditional.\\
\hline
\end{tabular}

\

\phantomsection
\label{GPY}
In fact, Goldston-Pintz-Y{\i}ld{\i}r{\i}m (GPY) were able to push their method even further to show unconditionally \cite{GoldstonPintzYildirim2010} that
$$
\liminf_{n\to\infty} {p_{n+1}-p_{n}\over (\log p_{n})^{1/2}(\log\log p_{n})^{2}}\ <\  \infty.
$$
So consecutive primes  infinitely often differ by  about square-root of the average gap. Moreover, assuming the primes have {\it any} level of distribution $\vt>1/2$, that is, any level in the \hyperref[conj:EH]{E-H} range (going beyond \hyperref[thm:BV]{B-V}), 
the GPY method gives a conditional proof of the Bounded Gaps Conjecture. 

The GPY method has been explained in great detail in a number of beautiful expositions (e.g. \cite{Sound2007, GoldstonPintzYildirim2007}) so we will not repeat the discussion here, contenting ourselves with just a few words on Zhang's advances. 
Once GPY was understood by the community, the big open question, in light of \hyperref[thm:BFI]{BFI}, was whether the ``weights'' $\gl(q)$ from \eqref{eq:BFI} could somehow be incorporated into the GPY method, so that in the resulting error analysis, \hyperref[thm:BV]{B-V} could be replaced by BFI.
There was even a meeting at the American Institute of Mathematics in November 2005, at which one working group was devoted to exactly this problem. At the time, at least to 
some, 
it did not seem promising.

Yitang Zhang's accomplishment, then, was threefold. He first changed the GPY weighting functions in a clever way  (in fact a similar change had been observed independently by Motohashi-Pintz
\cite{MotohashiPintz2008} and others),  then he proved an analogue of  the GPY sieving theorem with his new weights (as Motohashi-Pintz had also done), and finally (and most spectacularly!), he proved a more flexible\footnote{The most important aspect of Zhang's version is that the shift variable $a$ is allowed to vary, as opposed to \eqref{eq:BFI} where it must be fixed.}   
analogue  of \hyperref[thm:BFI]{BFI} 
which incorporates 
his
 new
 weights.
 In this technical tour-de-force, he was able to break the $\vt=1/2$ barrier in the weighted level of distribution of the primes, and complete the program initiated by Goldston in \cite{Goldston1992}.

Here is one final Basic Notion on this topic: \hyperref[thm:Zhang]{Zhang's Theorem}, at least as it currently stands, is ineffective! What he actually proves is a twin-prime analogue of Bertrand's Postulate (that for any $x>1$, there is a prime between $x$ and $2x$).

\vskip.06in
\noindent
{\bf Zhang's Theorem, Again:}
{\it 
For every $x$ sufficiently large, there is a pair of primes with difference at most $7\times 10^{7}$ in the range $[x,2x]$.}

\vskip.06in
How large is sufficiently large? It depends on whether or not GRH is true! Like most others, Zhang too relies at some early stage on the ineffective \hyperref[thm:SiegWalf]{Siegel-Walfisz Theorem}, and for this reason cannot 
escape 
Siegel zeros. 
(On the other hand, Heath-Brown \cite{HeathBrown1983a} has famously shown that if GRH fails and there {\it is} a particularly ``bad'' sequence of  Siegel zeros,
 then the 
Twin Prime Conjecture would follow!)

For further reading, we recommend any number of excellent texts, e.g. \cite{DavenportBook, IwaniecKowalski, FriedlanderIwaniecBook}, and of course the original papers.

\

{\bf Added in proof:}
It is very fortunate for the  author that he chose to 
focus
this survey 
 on the ``level of distribution'' aspect of Zhang's work. In November 2013, James Maynard \cite{Maynard2013}, developing an earlier attempted version of a method by Goldston-Y\i ld\i r\i m, succeeded in proving the even more shocking result:
 
\vskip.06in
\noindent
{\bf Maynard's Theorem:}
{\it 
For any $\ell\ge1$,}
$$
\liminf_{n\to\infty}(p_{n+\ell}-p_{n})<\infty.
$$

\vskip.06in
That is, one can find not only prime pairs which differ by a bounded amount, but also prime triples, quadruples, etc.  
Most remarkably, Maynard only needs the primes to have {\it any} exponent of distribution $\vt>0$ for his method to work (so now not even 
\hyperref[thm:BV]{B-V} is needed)!
Nevertheless, Zhang's spectacular achievement in going beyond the Riemann hypothesis in giving a flexible (weighted) exponent of distribution beyond $\gt=1/2$
will stand the test of time, and will surely find other applications. For a beautiful exposition of this aspect of Zhang's work, see the recent arXiv posting by Friedlander-Iwaniec \cite{FriedlanderIwaniec2014}. 

\newpage

\section{The Affine Sieve}

The goal of the Affine Sieve, initiated by Bourgain-Gamburd-Sarnak and completed by Salehi Golsefidy-Sarnak, is to extend to the greatest generality 
possible the mechanism of the Brun sieve. We will first discuss the former in \S\ref{sec:Brun} before turning our attention to the latter. In \S\ref{sec:Area} we motivate the general theory with an elementary problem, before presenting (some aspects of) the general theory in \S\ref{sec:genAff}.
This will again not be a 
rigorous or
comprehensive survey (for which  we refer the reader to any number of expositions, e.g. \cite{Sarnak2008a, SarnakNotes2010, Green2010, Kowalski2011,  Salehi2012}, in addition to the original papers \cite{BourgainGamburdSarnak2006, BourgainGamburdSarnak2010,  SalehiSarnak2011}), but rather  (we hope) a gentle introduction for the beginner. 
We apply the general theory to a few more illustrative examples in \S\S\ref{sec:genApply}--\ref{sec:Thin},
where we also give a discussion of Thin Orbits, see in particular \S\ref{sec:thin}.

While
the general theory is in principle ``complete,'' in that the Brun sieve can now be executed on matrix orbits, the whole program is far from 
finished, if one wishes to produce actual primes or almost-primes with very few factors in specific settings. 
We wish to highlight here some 
instances
 in which
 one can go beyond the capabilities of the general theory.
In certain special settings, one can now produce actual primes in Affine Sieve-type problems by applying
a variety of methods, each
 completely different  
from the general framework; we review some of these in
%
\S\ref{sec:Ex}.
Finally, we
discuss
in \S\ref{sec:AffLevels}
  other special
settings in which,
though
primes cannot yet be produced,
novel techniques 
have
nevertheless
given 
improved
levels of distribution, in the end coming quite close to producing primes.
We hope these give the reader some sense of the present landscape.


\subsection{The Brun Sieve}\label{sec:Brun}\

As throughout, we give only the most basic ideas.
The first sieving procedure for producing tables of primes is credited to the ancient  Eratosthenes ($\sim$200 BCE), whose
method exposes
  a simple but important observation: if $n< x$ and $n$ has no prime factors below $\sqrt x$, then $n$ is prime. Thus to make a table of the primes up to 100, one needs only to strike out (sieve) numbers divisible by 2, 3, 5, and 7 (the primes below $10=\sqrt{100}$).
A very slight generalization of the above is that: if $n<x$ has no prime factors below $x^{1/(R+1)}$, then $n$ is a product of at most $R$ primes. We call such 
a number $R$-almost-prime,
 and let $\cP_{R}$ denote the set of $R$-almost-primes. 

As a warmup, let us try (and fail) to prove the \hyperref[thm:PNT]{PNT} by sieving. 
To count the primes up to $x$, we first take the integers up to $x$ (there are $x$ of them), throw out those divisible by $2$ (there are {\it roughly} $x/2$ of them), then $3$ (there are roughly $x/3$ of them), and so on for all primes up to $\sqrt x$. But then we have twice thrown out multiplies of $2\times3=6$, so should add them back in (there are roughly $x/6$ of them), and so on goes the familiar inclusion-exclusion principle:
$$
\pi(x)\ \overset{?}{\approx}\  x-{x\over 2}-{x\over 3}-\cdots+{x\over 2\times3}+{x\over 2\times5}+\cdots-{x\over2\times3\times5}-\cdots
$$ 

The problem with this approach is two-fold. First of all, the word ``roughly'' above is very dangerous; it hides the 
error
$r_{q}=r_{q}(x)$ in 
\be\label{eq:nUpXmodQ}
\#\{n<x:n\equiv0(\mod q)\}\ = \
{x\over q}+r_{q},\qquad\qquad 
|r_{q}|<1.
\ee
These remainder terms $r_{q}$, when added together with absolute values in the inclusion-exclusion procedure, very quickly swamp the main term.
Perhaps we are simply too crude and a better estimation of these errors can make the above rigorous? Alas, were this the case, an elementary analysis (see, e.g., \cite{Granville1995}) will predict 
that
$$
\pi(x)\ \sim \ 2e^{-\g}{x\over \log x},
$$
where $\g$ is the Euler-Mascheroni constant, 
$$
\g\ :=\ \lim_{n\to\infty}\left(1+\frac12+\frac13+\frac14+\cdots+\frac1n - \log n\right)\ \approx\ 0.577.
$$
Since $2e^{-\g}\approx1.12$, we would be off by a constant from the truth. So these error terms {\it must} be at least of the same order as the main term. This simple sieving procedure cannot by itself prove the PNT. 

It was a technical tour-de-force when  Brun \cite{Brun1919} managed to push arguments of the above flavor, together with a heavy dose of combinatorics and other estimates, to prove 
a different
type of approximation to the Twin Prime Conjecture, in some sense orthogonal to \hyperref[thm:Zhang]{Zhang's Theorem}.

\vskip.1in
\noindent
{\bf Brun's Theorem (1919):}
{\it There are infinitely many integers $n$ so that both $n$ and $n+2$ have at most nine prime factors.}

\vskip.1in
That is, infinitely often $n$ and $n+2$ are simultaneously in $\cP_{R}$ with $R=9$.
After much work by many people, the sieve was finally pushed to its limit\footnote{See the discussion of the ``parity problem'' in, e.g., \cite{FriedlanderIwaniec2009a}.} in \cite{Chen1973}:

\phantomsection
\vskip.1in
\noindent
{\bf Chen's Theorem (1973):}\label{thm:Chen}
{\it There are infinitely many primes $p$ so that $p+2$ is either prime or the product of two primes.}

\vskip.1in
It was realized long ago that this 
sieving procedure
 applies to much more general problems.
%
Suppose we have an infinite set of %
natural numbers 
\be\label{eq:setcS}
\cS\subset
\N
\ee 
and 
wish
 to prove the existence and abundance of  primes or $R$-almost-primes in $\cS$.
Roughly speaking, 
 all that is needed is
 an appropriate analogue of \eqref{eq:nUpXmodQ}. 
In particular, suppose
 that
$\cS$ is fairly 
well
 distributed on average among multiples of $q$, in the sense that
 \footnote{There are various ways the assumption \eqref{eq:cSmodq} can (and should) be relaxed and generalized further, but for the purposes of our discussion, we will ignore all technicalities and stick with this simple-minded version.}
\be \label{eq:cSmodq}
\#\{n\in\cS\cap[1,x]:n\equiv0(q)\}
\ =\
\frac1q
{\#\{\cS\cap[1,x]\} 
}
+
r_{q}
,
\ee
(or perhaps with $1/q$ in \eqref{eq:cSmodq} replaced by 
some analytically similar function like $1/\phi(q)$
),
where
the errors 
are
controlled by
 \be\label{eq:rqErr}
\sum_{q<Q}
\left|
r_{q}
\right|
\ =\
o\left(
\#\{\cS\cap[1,x]\}
\right)
,
 \ee
for some $Q$.
Such an expression should look familiar; it is in some sense the generalization of \eqref{eq:LevelDef}, and $Q$ is likewise called a {\it level of distribution} for $\cS$. 
Then the sieve technology (again 
very
roughly) tells us that if $Q$ can be taken as large as a power of $x$, say 
\be\label{eq:Qsieve}
Q\ = \ x^{\vt-\vep}
\ee
for some {\it exponent of distribution} $\vt>0$, then 
$\cS$ contains $R$-almost-primes, with 
\be\label{eq:Rsieve}
R
\ = \
\left\lceil\frac1\vt+\vep\right\rceil.
\ee 

For example, if $\cS$ is the set of shifted primes, $\cS=\{p+2:p\text{ prime}\}$, then the \hyperref[thm:BV]{Bombieri-Vinogradov Theorem} gives 
us an exponent of distribution $\vt=1/2$, which 
gives $R$-almost-primes in $\cS$ with $R
=
\lceil
2+\vep\rceil=3$. 
To obtain \hyperref[thm:Chen]{Chen's Theorem} 
is much much harder. 

\subsection{Affine Sieve Warmup: Pythagorean Areas}\label{sec:Area}
\

Arguably the oldest ``Affine Sieve'' problem is the following.
 Let $(x,y,z)$ be a Pythagorean triple, that is, an integer solution to the equation $x^{2}+y^{2}=z^{2}$.
What can one say about the number of prime factors of the 
area 
$
\foh xy$ 
of a Pythagorean triple?

It
was known
to the ancients that  Pythagorean triples $\bx=(x,y,z)$ with coprime entries and $x$ odd are parametrized by  coprime pairs $(c,d)$ of opposite parity with
\be\label{eq:PythParam}
x=c^{2}-d^{2},\qquad y=2cd,\qquad z=c^{2}+d^{2}.
\ee
In fact, it is easy to see that the area is always divisible by $6$, so we can 
further remove unwanted prime factors by 
studying
 the function 
$f(\bx)=\frac1{12}xy$.
Observe that
in the parametrization \eqref{eq:PythParam}, we have 
\be\label{eq:area}
f(\bx)=\frac1{12} xy=\frac16cd(c+d)(c-d).
\ee
When does $f(\bx)$ have 
few prime factors? That is,
for which $R$ and triples $\bx$ is $f(\bx)\in\cP_{R}$?

One can easily check that 
 there are only finitely many pairs $(c,d)$ so that 
  \eqref{eq:area} is 
 the product of two primes. 
The largest such pair is $(c,d)=(7,6)$, which corresponds to the triple $\bx=(13,84,85)$ of one-sixth area $f(\bx)=
\frac1{12}13\times84=
91=7\times13$.

\begin{figure}
\includegraphics[height=1.5in]{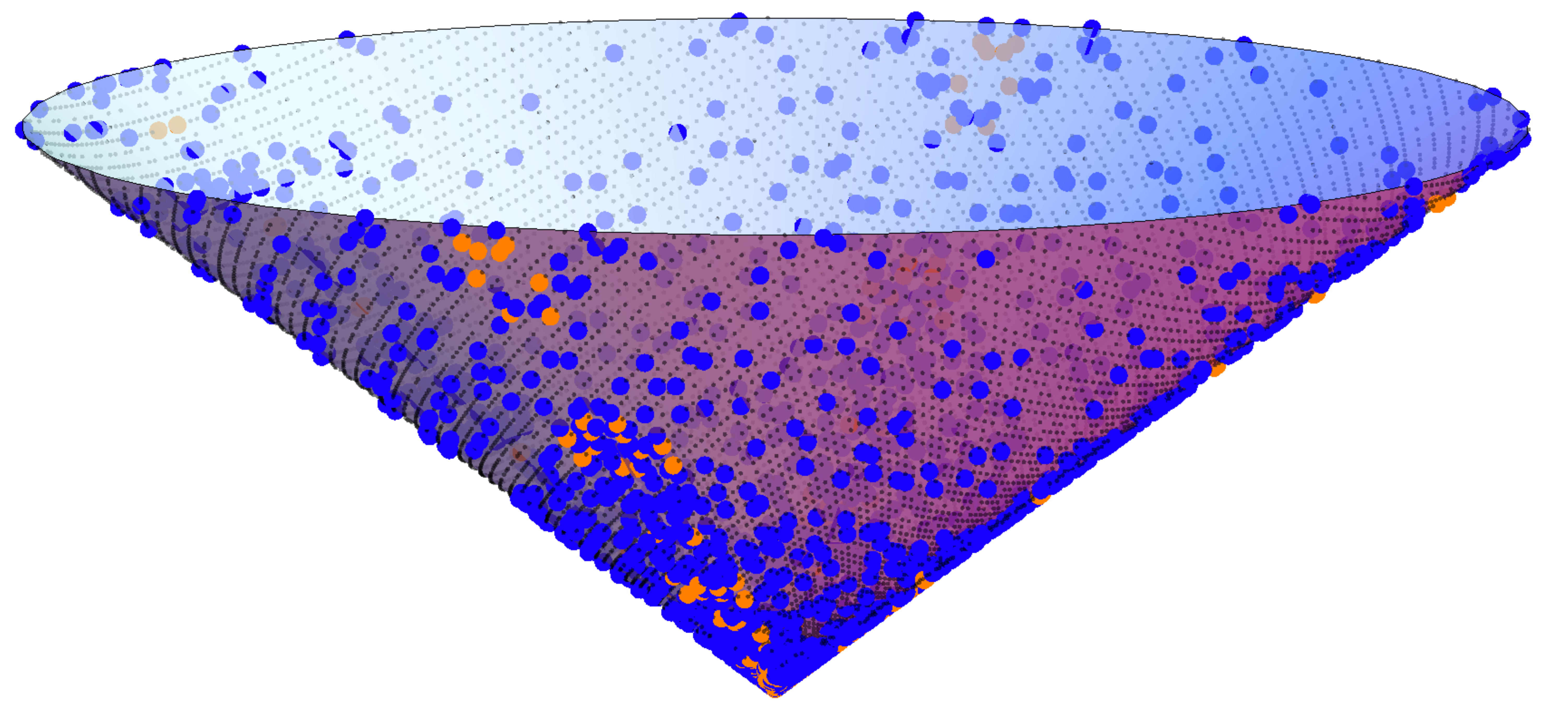}
\caption{A piece of the cone $V$ in \eqref{eq:VPyth} with markings at the 
  primitive Pythagorean triples $\bx$%
. Points $\bx
$ are marked according to whether the ``area''
$f(\bx)=\frac1{12}xy$
 is 
in $\cP_{R}$ with $R\le 3$ (\protect\includegraphics[width=.1in]{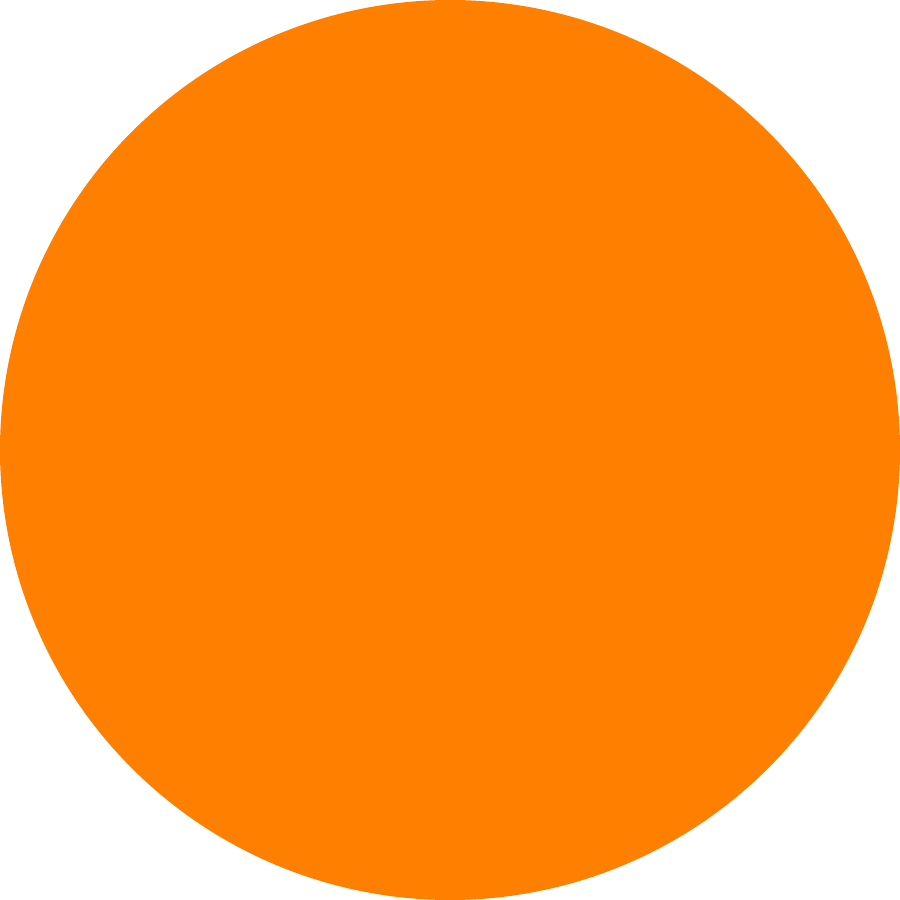}),
$R=4$ (\protect\includegraphics[width=.1in]{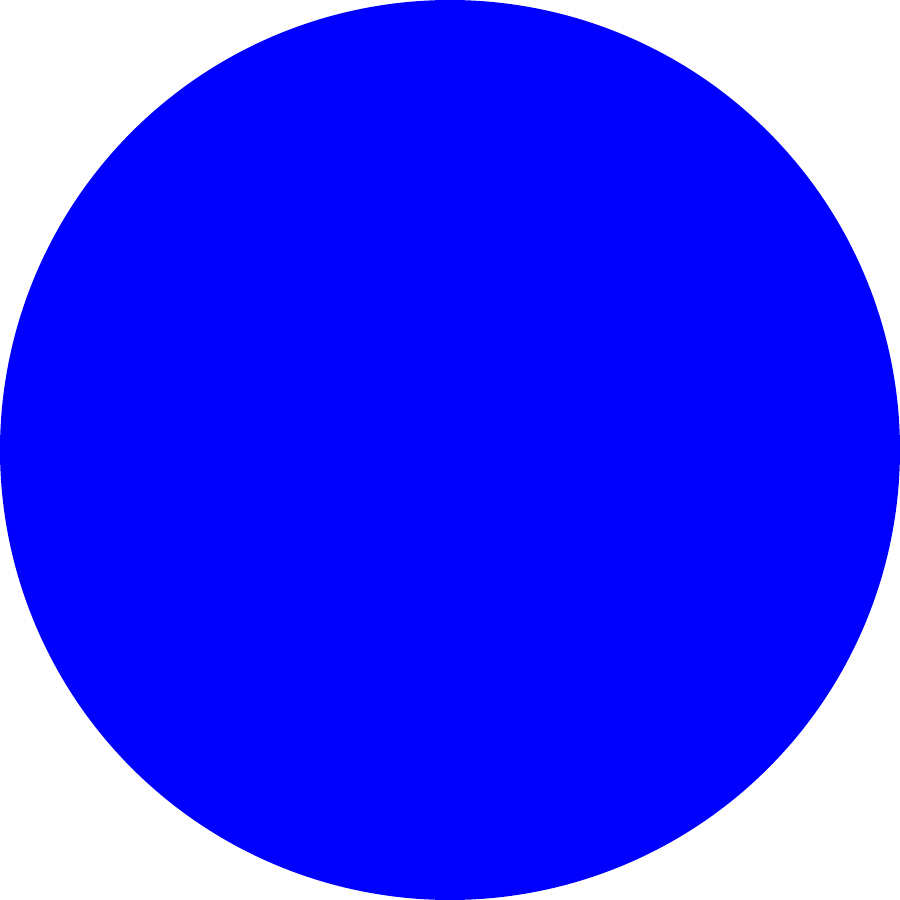}),
 or 
$R\ge5$ (
\protect\includegraphics[width=.1in]{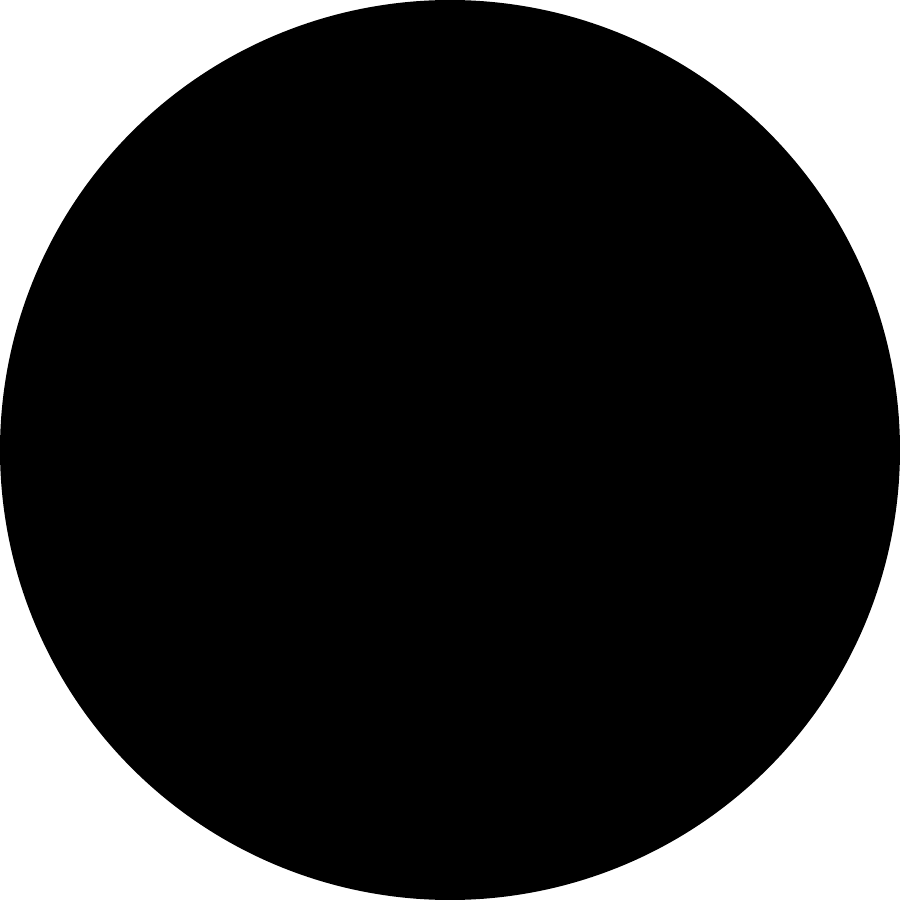})
  .}
\label{fig:PythA}
\end{figure}

Allowing $R=3$ primes
, we could set, say, $d=2$; then from  \eqref{eq:area}, we are asking for many $c$'s so that $\frac13c(c-2)(c+2)$ is 
  the product of three primes (a type of ``triplet prime'' problem). 
As the reader may surmise, it is expected that infinitely many such $c$'s exist, but this seems far outside the range of what can be proved today. Nevertheless, it should be clear that for $f(\bx)$ to have 
three prime factors, $\bx$ must be of some ``special'' form, so either $c$ or $d$ (or their sum or difference) must be ``small'';
see Figure \ref{fig:PythA}.

A better way of saying this is to use the Zariski topology. 
The ambient variety 
on which all Pythagorean triples live is the cone $V$ given  by 
\be\label{eq:VPyth}
V\ : \ F(\bx)=0,
\ee
where $F$ is the 
 quadratic form 
 \be\label{eq:FPyth}
F(\bx)\ =\ x^{2}+y^{2}-z^{2}
. 
\ee
Let $\cX_{R}$ denote the set 
of integer Pythagorean triples  $\bx$ with $f(\bx)\in\cP_{R}$.
Then a 
restatement of 
the ``smallness'' of points in $\cX_{3}$ is that
the (affine) Zariski closure of $\cX_{3}$ 
is a proper subvariety of $V$. That is, points in $\cX_{3}$ have extra algebraic relations. 

This ``smallness'' somehow fundamentally changes the nature of the problem; e.g. 
setting
$d=2$ as above, one is asking for a ``triplet prime'' type statement, rather than the original area problem. 
It seems natural, then,  to exclude such small solutions. 
That is, we shall insist on finding an $R$ so that the set $\cX_{R}$ is Zariski dense in $V$; this means that any polynomial which vanishes on all of $\cX_{R}$ must also vanish on $V$.

If we now allow $R=4$
 prime factors, 
then we see in Figure \ref{fig:PythA} that such points seem to spread out all over the cone. In fact, it was observed in \cite{BourgainGamburdSarnak2010} that Green-Tao's revolutionary work \cite{GreenTao2010} on linear equations in primes 
rigorously establishes
the Zariski density of $\cX_{4}$  in $V$. This is because  $f(\bx)$ in \eqref{eq:area} is the product of four linear factors in two variables, which in the Green-Tao nomenclature is a system of ``finite complexity'' (we refer the reader to their paper for the definition, which is not needed here).
Thus the problem of Pythagorean areas, at least if one insists on Zariski density, is completely solved.

\subsubsection{Reformulation}
\

What does this simple problem have to do with orbits? 
Let 
$$
G=\SO_{F}(\R)=\SO_{2,1}(\R)
$$ 
be the real special orthogonal group preserving
the quadratic form
 $F$ in \eqref{eq:FPyth}; that is,
\be\label{eq:GPyth}
G\ =\ \{g\in\SL_{3}(\R):F(g\cdot\bx)=F(\bx),\ \forall\bx\}.
\ee
This is a nice algebraic  (defined by polynomial equations
) Lie group, and its integer subgroup 
\be\label{eq:GamPyth}
\G\ :=\ \SO_{F}(\Z)
\ee 
is a nice arithmetic (the set of integer points on an algebraic group) discrete group.

To make these groups slightly less mysterious, it is a well-known fact (see, e.g., the discussion in \cite[\S4]{Kontorovich2013}) 
that they
 can be parametrized, as follows. It can be checked that whenever 
 $$
 \mattwo abcd\in\SL_{2}(\R),
 $$ 
the matrix
\be\label{eq:morph}
g:=
\left(
\begin{array}{ccc}
 \frac{1}{2} \left(a^2-b^2-c^2+d^2\right) & a c-b d & \frac{1}{2}
   \left(a^2-b^2+c^2-d^2\right) \\
 a b-c d & b c+a d & a b+c d \\
 \frac{1}{2} \left(a^2+b^2-c^2-d^2\right) & a c+b d & \frac{1}{2}
   \left(a^2+b^2+c^2+d^2\right)
\end{array}
\right)
\ee
is in $G$.
Likewise, $\G$ is essentially the image under the above morphism of the more familiar discrete group $\SL_{2}(\Z)$. 

The set of all primitive  Pythagorean triples (up to symmetry) is then given by the orbit:
\be\label{eq:orbit}
\cO\ :=\ \G\cdot\bx_{0},
\ee
where $\bx_{0}$ is any primitive
base point $\bx_{0}\in V(\Z),$ e.g. 
$$
\bx_{0}=(3,4,5).
$$
To study the 
area, we 
again 
consider the function
 $f(\bx)=\frac1{12}xy$, and ask for an $R$ so that
 the set $\cX_{R}$ of $\bx\in\cO$ with $f(\bx)\in\cP_{R}$ 
 is Zariski dense in
 the cone $V$ in \eqref{eq:VPyth}, which is
  the Zariski closure of $\cO$. 

\subsection{
The General Procedure}\label{sec:genAff}\

We have taken a very simple problem and made it look very complicated. 
But now we have seen almost all of the essential features of the general 

\vskip.1in
\noindent
{\bf Affine Sieve:} 
One takes
\begin{enumerate}
 \item  a finitely generated subgroup $\G$  
 of $\GL_{n}(\Q)$
(later we will want to relax this to allow 
 semigroups), 
 \item 
  some base point $\bx_{0}\in\Q^{n}$, which then forms the orbit $\cO$ as in \eqref{eq:orbit}, and
\phantomsection
  \item\label{(3)}
a polynomial function $f$ which takes integer values on $\cO$. 
  \end{enumerate}
With this data, one asks for an (or the smallest) integer $R<\infty$ so that
the set 
$$
\cX_{R}=\cX_{R}(\cO,f)
\ :=\
\{\bx\in\cO:f(\bx)\in\cP_{R}\}
$$ 
is Zariski dense in the Zariski closure of $\cO$.
In practice, the Zariski density is not hard to establish, so we will simply say that $f(\cO)$ 
contains $R$-almost-primes (or that we have produced $R$-almost-primes) to mean the more precise 
statement.

\vskip.1in
Let us see now how the general Affine Sieve method proceeds.
In the notation of \eqref{eq:setcS}, 
we wish to sift for $R$-almost-primes
in the set
$$
\cS\ :=\ f(\cO).
$$
As in \eqref{eq:cSmodq}, we must understand the distribution of $\cS\cap[1,x]$ among the multiples of $q$ up to some level $Q$. 
Roughly speaking, if $\g\in\G$  is of size $\|\g\|$ 
about
$T$, then so is the size of $\|\bx\|$, where $\bx=\g\cdot\bx_{0}$, 
since the base point $\bx_{0}$ is fixed.%
\footnote{%
For simplicity, take all norms here to be Euclidean,
though
in many settings it is advantageous
(or even necessary, since we do not yet know how to count with archimedean norms in full generality!)
 to use other norms, 
e.g.,
the wordlength metric in the generators of $\G$.}
 If $f$ is a polynomial of degree $d$, then generically $f(\bx)$ is of size $T^{d}$ for such an $\bx$. Hence restricting $\cS$ to $f(\bx)<x$ is roughly the same as restricting $\|\g\|<T$ with $T= 
x^{1/d}$. 
The left hand side of \eqref{eq:cSmodq} 
may
then 
be
 captured in essence  by
 \be\label{eq:1003s}
\sum_{\g\in\G\atop\|\g\|< x^{1/d}}\bo_{\{f(\g\cdot\bx_{0})\equiv0(q)\}}.
 \ee
We should first determine what 
happens if
$q=1$, that is, when the congruence condition is dropped. Say the group $\G$ has {\it exponent of growth} 
\be\label{eq:gdIs}
\gd\ > \ 0,
\ee 
which means roughly that the number of points in $\G$ of norm at most $T$ is about $T^{\gd}$, or
 \be\label{eq:1003ss}
\sum_{\g\in\G\atop\|\g\|< x^{1/d}}1
\ 
=
 \
x^{\gd/d+o(1)}.
\ee
Note that
for general $q$, the condition $f(\g\cdot\bx_{0})\equiv0(q)$ is only a restriction on $\g$ mod $q$, so we can decompose the sum above into residue classes as
 \be\label{eq:1003sss}
\sum_{\g_{0}\in\G(\mod q)}
\bo_{\{f(\g_{0}\cdot\bx_{0})\equiv0(q)\}}
\left[
\sum_{\g\in\G\atop\|\g\|< x^{1/d}}\bo_{\{\g\equiv\g_{0}(q)\}}
\right]
.
\ee
The bracketed term above is the key to the whole game. 
What do we expect? 
If the euclidean ball  in $\G$ of size $x^{1/d}$ is equidistributed among the possible residue classes mod $q$, then the bracketed term should be ``roughly'' equal to
$$
\frac1{|\G(\mod q)|}
\sum_{\g\in\G\atop\|\g\|< x^{1/d}}1
.
$$
In reality, 
one can prove today in 
some
generality thanks to the work of many people (e.g. \cite{
Selberg1965, 
LaxPhillips1982,
LuoRudnickSarnak1995,
Kim2003,
BlomerBrumley2011,
BurgerSarnak1991, Clozel2003, 
SarnakXue1991, 
Gamburd2002,
BourgainGamburd2008,
BourgainGamburdSarnak2010,
BourgainGamburdSarnak2011,
Helfgott2008,
BreuillardGreenTao2011,
PyberSzabo2010,
SalehiVarju2012})
an estimate of the form
 \be\label{eq:1003ssss}
\sum_{\g\in\G\atop\|\g\|< x^{1/d}}\bo_{\{\g\equiv\g_{0}(q)\}}
\ =\ 
\frac1{|\G(\mod q)|}
\sum_{\g\in\G\atop\|\g\|< x^{1/d}}1
+
O\left(q^{C}
\left[
\sum_{\g\in\G\atop\|\g\|< x^{1/d}}1
\right]^{1-\gT}
\right)
.
\ee
Here $C<\infty$ and $\gT\ge0$ are some constants, and if
\be\label{eq:gTis}
\gT\ >\ 0,
\ee
then $\gT$ is
 often referred to as a ``spectral gap'' for $\G$.
%
%
Results of this type follow (with quite a bit of work in many separate cases) from theorems (or partial results towards conjectures)
going under various guises; 
some of these ``buzzwords'' are:
the Selberg $1/4$-Conjecture, the generalized Ramanujan conjectures, 
mixing rates for homogeneous flows, temperedness of representations, 
resonance-free regions for transfer operators, 
expander graphs, among many others; see, e.g. \cite{Sarnak1995, Sarnak2004, Sarnak2005,  HooryLinialWigderson2006, Lubotzky2012, BlomerBrumley2013}.

Now inserting \eqref{eq:1003ssss} into \eqref{eq:1003sss},  using \eqref{eq:1003ss}, and assuming that the proportion of $\g_{0}$ in $\G(\mod q)$ with 
$$
f(\g_{0}\cdot\bx_{0})\equiv0(q)
$$ is about $1/q$ (for example, $f$ should
 not be identically zero%
 ), we obtain an estimate for \eqref{eq:1003s} roughly of the form \eqref{eq:cSmodq}, with
\be\label{eq:rqBnd}
|r_{q}|\ \ll\  q^{C}x^{\gd(1-\gT)/d}.
\ee
(The value of the constant $C$ may change from line to line.) Again using \eqref{eq:1003ss} as an approximation for $\#\cS\cap[1,x]$, we 
obtain
 that \eqref{eq:rqErr} holds with
$$
Q\ = \ X^{\gd\gT/(Cd)-\vep},
$$
say, for any $\vep>0$. 
Thus the set $\cS$ has exponent of distribution 
\be\label{eq:expDist}
\vt \ = \ {\gd\gT\over Cd},
\ee
and hence contains 
$R$-almost-primes 
with
\be\label{eq:Rgen}
R\ =\ \left\lceil {Cd\over \gd\gT}+\vep\right\rceil.
\ee
So
as long as $C<\infty$, that is, the dependence on $q$ in the error term of \eqref{eq:1003ssss} is at worst polynomial, and as long as the ``spectral gap''  $\gT$ is strictly 
positive, this general sieving procedure produces $R$-almost-prime values in $\cS$ for some $R<\infty$.

\pagebreak


\subsection{Applying the General Procedure}\label{sec:genApply}

\subsubsection{Fibonacci Composites}\

Again it is instructive to first see how the general method can fail to work. 
Let $\G$ be the semigroup generated by the square of the matrix $\mattwos 0111$, set $\bx_{0}=(0,1)$, with orbit $\cO=\G\cdot\bx_{0}$, and consider the function $f(x,y)=x$. It is elementary to check that here 
$$
\cS
\ =\ 
f(\cO)
\
=
\
\{{\ff}_{2n}\}
$$ 
is just the set of even-indexed Fibonacci numbers, $\ff_{2n}$. This 
set
 is much too thin for the above methods to apply, since the number of Fibonacci numbers up to $x$ is about $\log x$; that is, the group $\G$ has exponent of growth $\gd$ in \eqref{eq:gdIs} equal to zero.
In particular, a counting result of the type \eqref{eq:1003ssss} is simply impossible, and one cannot establish a positive exponent of distribution $\vt$ as in \eqref{eq:expDist}. 

In fact, there seems to be good reason for the sieve to fail in this context. While it is believed that infinitely many Fibonacci numbers are prime, there is some heuristic evidence that if $n$ is composite, then the $n$th Fibonacci number $\ff_{n}$ has at least on the order of $n$ prime factors. 
Assuming this heuristic, there should not exist a finite $R$ for this setting; that is, the sieve does not work here because it must not. 
(Note that the Zariski closure of the group $\G$ here is a torus, $\C^{\times}$; as Sarnak likes to say, for the Affine Sieve, ``the torus is the enemy!'')

\subsubsection{Back to Pythagorean Areas}\label{sec:PythA}

What does the above 
procedure give for Pythagorean areas? The function $f(\bx)=\frac{1}{12} xy$ is quadratic, so $d=2$. It is not hard to see that $\G$ has growth exponent $\gd$ in \eqref{eq:gdIs} equal to $1$. Selberg's $1/4$-Conjecture, if true, would  imply an estimate (in smooth form) for \eqref{eq:1003ssss} with
 ``spectral gap''
  $\gT=1/2$; this is again 
a
square-root cancellation type phenomenon. 
Unconditionally, the best-known bound (due to Kim-Sarnak) 
proves \eqref{eq:1003ssss} with
$\gT=\frac12-\frac7{64}.$
The value for $C$ coming from  (a slight variant of) the above procedure can be whittled down to $2$. 
One small technicality is that our $f$  in \eqref{eq:area} is now the product of four irreducible factors, so  the fraction $1/q$ on the right hand side of \eqref{eq:cSmodq} should be replaced by $4/q$ (giving a sieve of ``dimension'' $4$); the sieve still works in the same way, just with a worse dependence of $R$ in \eqref{eq:Rsieve} on the level of distribution in \eqref{eq:Qsieve}. 
%
%

The above technicalities aside,
all this machinery 
will in the end produce 
an exponent of distribution $\vt$ of about $
1/10
$, and about 
\be\label{eq:Rarea}
R\ =\ 30
\ee 
primes,
%
falling far short of Green-Tao's optimal result $R=4$.
Of course the orbit $\cO$ here is very simply described, making its study amenable to other means.
In the  following subsection, 
we give a sampling of problems in which more elementary descriptions do not seem advantageous (or even possible), yet where the Affine Sieve applies as just indicated.
We hope these
 serve to illustrate
 some of the
 power and robustness of the Affine Sieve.

\subsection{More Examples: Anisotropic and Thin Orbits}\label{sec:Thin}

\subsubsection{Anisotropic ``Areas''}\label{sec:aniso}\

Keeping a nearly identical setup, let us change ever so slightly the quadratic form $F$ from \eqref{eq:FPyth} to
$$
F(\bx)\ =\ x^{2}+y^{2}-3z^{2}
.
$$
The salient features of this form 
are that, like \eqref{eq:FPyth},
 it
 is
rational (the ratios of its coefficients are in $\Q$) and indefinite (it takes positive and negative values),
but
unlike \eqref{eq:FPyth}, it is {\it anisotropic} over $\Q$. This means that
%
it has no non-zero rational points on the cone $F=0$. (Exercise.) So to have an integral orbit, we can change our variety $V$ from \eqref{eq:VPyth} to, say,
$$
V\ : \ F(\bx)=1
,
$$
which over $\R$ is a one-sheeted hyperboloid containing the integer base point $\bx_{0}=(1,0,0)$. Let 
$G=\SO_{F}(\R)$ 
now
be the real special orthogonal group preserving this new form, and let $\G=\SO_{F}(\Z)$ be the arithmetic group of integer matrices in $G$. Taking 
the orbit $\cO=\G\cdot\bx_{0}$ and function 
$$
f(\bx)=\foh xy
$$ 
as an analogue of ``area,'' one can compute (see \cite{Kontorovich2011}) that $\cS=f(\cO)$ is essentially the set of all values of 
\be\label{eq:anisArea}
(a^{2}-b^{2}+3c^{2}-3d^{2})
(ab+3cd),
\ee
where $a,b,c,d$ range over all integers satisfying
\be\label{eq:anis}
a^{2}+b^{2}-3c^{2}-3d^{2} \ = \ 1.
\ee
(In fancier language, the spin group of $\G$ is isomorphic to the norm one elements of a particular quaternion division algebra.)

Needless to say, the Green-Tao technology of linear equations 
is not designed to
handle this new set $\cS$, while the Affine Sieve works in exactly the same way as previously described (in this setting, it was executed by Liu-Sarnak \cite{LiuSarnak2010}), producing $R$-almost-primes\footnote{%
For the experts, this number is about half of that in \eqref{eq:Rarea}, due to \eqref{eq:anisArea} being a two-dimensional sieve problem instead of \eqref{eq:area} which is four-dimensional. In the anisotropic case, there are no ``extra'' parametrizations
like \eqref{eq:PythParam}, so the ``area'' is only a product of two irreducible factors, not four.}
 with $R=16$. 

\subsubsection{A Thin Group}\label{sec:thin}
\

While the 
group $\G$ 
in \S\ref{sec:aniso}
was more complicated, it was still arithmetic; 
in particular,
any solution in the integers to the polynomial equation \eqref{eq:anis} gave (by 
a simple formula) an element in $\G$. The situation is 
even
more 
delicate
if the group $\G$ is 
restricted
to 
some {\it infinite} index subgroup of $\SO_{F}(\Z)$. Here is a quintessential ``thin'' (see below for the definition) group.

\begin{figure}
\includegraphics[height=1.5in]{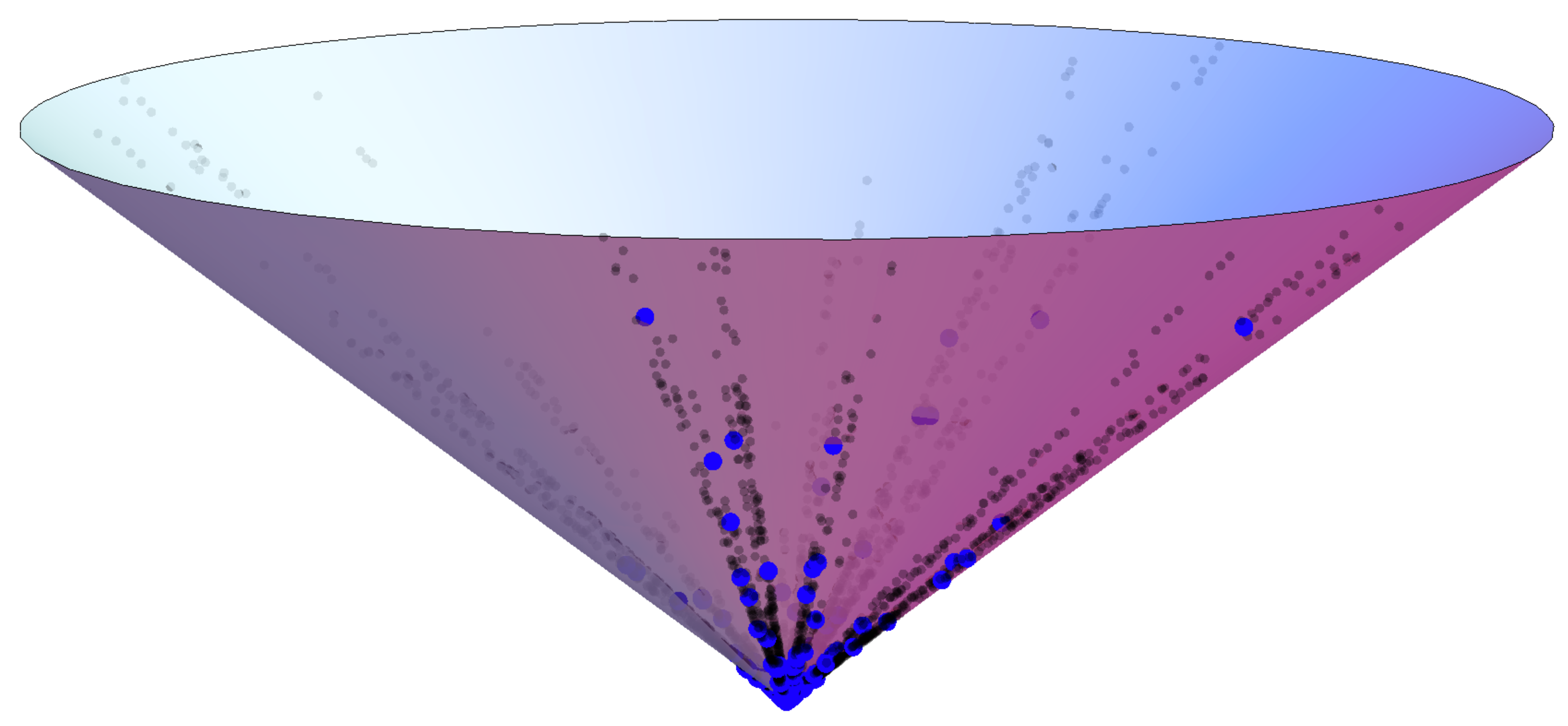}
\caption{A piece of the thin
Pythagorean 
orbit $\cO$ 
in \eqref{eq:cOthin}.
 Points $\bx\in\cO$  are again marked according to whether the ``area''
$f(\bx)=\frac1{12}xy$
 is  
in $\cP_{R}$ with $R\le 3$ (\protect\includegraphics[width=.1in]{DotOrange.pdf}),
$R=4$ (\protect\includegraphics[width=.1in]{DotBlue.pdf}),
 or
$R\ge5$ (\protect\includegraphics[width=.1in]{DotBlack.pdf}).
}
\label{fig:PythThin}
\end{figure}

Let us return again to the Pythagorean setting of \eqref{eq:FPyth} and the cone \eqref{eq:VPyth} with base point $\bx_{0}=(3,4,5)$ and the ``area'' function $f(\bx)$ in \eqref{eq:area}. 
For the sake of being explicit,
let $\G$ be the group generated by the two matrices
\be\label{eq:MsP}
M_{1}:=
\left(
\begin{array}{ccc}
 -7 & -4 & -8 \\
 4 & 1 & 4 \\
 8 & 4 & 9
\end{array}
\right)
\quad\text{and}\quad
M_{2}:=
\left(
\begin{array}{ccc}
 -1 & 2 & 2 \\
 -2 & 1 & 2 \\
 -2 & 2 & 3
\end{array}
\right)
,
\ee
which one can check are the images under the morphism \eqref{eq:morph} of 
$\mattwos1401$
and 
$\mattwos1021$, respectively. 
The orbit  
\be\label{eq:cOthin}
\cO=\G\cdot\bx_{0}
\ee
of $\bx_{0}$ under this 
group $\G$ is illustrated in Figure \ref{fig:PythThin}; this is the picture one may keep in mind when thinking of thin orbits.

Note that, unlike Figure \ref{fig:PythA},  there   now seem to be only finitely many $\bx\in\cO$ with $f(\bx)$ having $R\le3$ prime factors; these are invisible at the scale drawn in Figure \ref{fig:PythThin}.
The (presumably infinite number of) points of ``special'' form visible in Figure \ref{fig:PythA} seem to disappear for this thin orbit, again reinforcing our 
suggestion
that Zariski density is the ``right'' 
demand
for the general setting.
\\

What do we mean by ``thin''?
There are a number of competing definitions of this word,  and we will need to give a new one to suit our purposes. 
The meaning of thin typically involved in the Affine Sieve refers to ``thin matrix groups'' (not to be confused with ``thin sets,'' as defined by Serre \cite[\S3.1]{Serre2008}),
which are finitely generated groups $\G<\GL_{n}(\Z)$ which have infinite index in the group of integer points of their Zariski closure.
That is, let 
$$
G=\Zcl(\G)<\GL_{n}
$$ 
be the Zariski closure of $\G$, and $G(\Z)$ its integer points; then $\G$ is called a thin matrix group if the index 
$$
[G(\Z):\G]=\infty.
$$ 

For our purposes, we
will
want to allow $\G$ to be a finitely generated semi-group of $\GL_{n}(\Z)$, but not necessarily a group.
In this case, we cannot speak of index, and need a different condition to characterize what should be considered  thin.
%
Moreover, we will want
   the flexibility to apply the adjective thin to either the (semi-)group $\G$, or the resulting orbit $\cO$, or the resulting set of integers $\cS=F(\cO)$. 
Our characterization will simply be by an archimedean degeneracy in the algebro-geometric closure, as follows.
   


\vskip.1in
\noindent
{\bf Definition: Thin Integer Set.} 
Let $\cZ\subset\Z^{n}$ be a set of integer vectors, let 
$\Zcl(\cZ)$ be the Zariski closure of $\cZ$, and let $B_{x}$ be a ball of radius $x>0$ (with respect to any fixed archimedean norm) about the origin in $\R^{n}$. 
We will call $\cZ$ a {\it thin integer set} if 
$$
\#(\cZ\cap B_{x}) \ =\
o\Big(
\#(\Zcl(\cZ)\cap\Z^{n}\cap B_{x})\Big) 
,
\qquad\text{ as }
x\to\infty.
$$
That is,  $\cZ$ has zero ``density'' inside the integer points of its Zariski closure.

\vskip.1in
It is an easy fact%
\footnote{%
A sketch for the experts: the trivial representation does not weakly occur in the regular action of $G=\Zcl(\G)$
on $L^{2}(\G\bk G)$ if and only if $\vol(\G\bk G)=\infty$, in which case Howe-Moore gives the decay of matrix coefficients. On the other hand, the count for arithmetic groups is known already by methods of Duke-Rudnick-Sarnak and Eskin-McMullen.
}
 that when $\G<\GL_{n}(\Z)$ is  group, then it is a thin matrix group if and only if it is a thin integer set in 
$ \Z^{n\times n}\cong\Z^{n^{2}}$. (Thanks to Peter Sarnak for insisting that we make our definition so that the two definitions would agree on their intersection.)

Our
group $\G=\<M_{1},M_{2}\>$ from \eqref{eq:MsP} has infinite index in $\SO_{F}(\Z)$, 
so is thin. 
Its exponent of growth can be estimated as 
$$
\gd\ \approx \ 0.59
\cdots
,
$$ 
which, it turns out, is also the Hausdorff dimension of
 the {\it limit set} of $\G$. The latter is roughly speaking the Cantor-like fractal set seen at the boundary at infinity in Figure \ref{fig:PythThin}; that is, the set of directions in which the orbit $\cO$ grows.

Now there is certainly no hope of a 
more direct 
approach to studying $\cS=f(\cO)$, as we cannot even determine, given a matrix $M\in\SO_{F}(\Z)$, whether it
is in the group $\G$.
Unlike the arithmetic group case, it is not enough to check whether the entries of $M$ satisfy some polynomial equations; instead 
one must determine whether $M$ can be realized as some word in the generators \eqref{eq:MsP}. 
As the general membership problem in a group is undecidable \cite{Novikov1955}, we had better avoid this issue.
Luckily, the 
standard
Affine Sieve procedure works just as described in \S\ref{sec:genAff}. (In this setting, the details were worked out by the author \cite{MyThesis, Kontorovich2009}, and the author with Oh \cite{KontorovichOh2012}). 

A good question to ask at this point 
might be: Why would anyone care about these 
strange thin groups? 
Here are
 just
 two motivations:
(1)  
thin
groups are in some sense ``generic'' (see, e.g., \cite{FuchsMeiriSarnak2012, Fuchs2012, Sarnak2014}, 
for a  discussion 
into  which we will not delve here), 
 and (2) many naturally-arising and interesting  
 problems
  {\it require} their study. 
Let us postpone our discussion of these natural problems for a moment, turning now to another topic. 

\subsection{The Affine Sieve 
Captures Primes}\label{sec:Ex}\

We
have
 described  the general procedure and explained how it works in a number of sample settings, 
but it is clear that without further ingredients, producing primes seems hopeless. 
 %
Yet, as we have already seen
in the case of Pythagorean areas, the Green-Tao theorem, using completely different tools, goes far beyond the present capabilities of the Affine Sieve. We give here but a sampling of four more settings in which other technologies  prove more successful, producing a minimal number of prime factors.

\subsubsection{
Matrix 
 Ensembles
 with Prime Entries}\

Now that we 
appreciate
the utility of posing problems in
terms of matrix orbits, 
why not
ask the following even  simpler Affine Sieve-type question: 
Among the 
set of all $n\times n$ integer matrices of, say, fixed determinant $D\ge1$, are 
 there 
 infinitely many
  with
  all  entries  prime? For example, here is a prime $3\times 3$ matrix of determinant $D=4$:
\be\label{eq:soln3}
  \det
  \left(
\begin{array}{ccc}
 3 & 5 & 7 \\
 11 & 13 & 17 \\
 5 & 13 & 19
\end{array}
\right)
\ =\ 4.
\ee

  How is this an Affine Sieve problem?
Let $V_{n,D}(\Z)$ be the set 
in question of all $n\times n$ integer matrices of determinant $D$.
The full\footnote{We will sometimes use ``full'' 
as the negation of ``thin.''}  
group $\SL_{n}(\Z)$ acts on $V_{n,D}(\Z)$ on  the left (determinant is preserved),
and a theorem of
 Borel and Harish-Chandra
tells us 
that $V_{n,D}(\Z)$ breaks up into finitely many such orbits.
  Thus we may as well just take one fixed matrix $M_{0}\in V_{n,D}(\Z)\subset \Z^{n^{2}}$ and consider the orbit
$$
\cO\ =\ \SL_{n}(\Z) 
\cdot
M_{0}
.
$$
 For an $n\times n$ integer matrix $M=(m_{ij})\in\cO$, our function $f$ is now the product of all coordinates, 
 $$
 f(M)\ =\ \prod_{i,j} m_{ij},
 $$ 
 which, being a product of $n^{2}$ terms, we would like to make $R$-almost-prime with $R=n^{2}$. 
 \\
 
 First let us consider the case $n=2$, that is, for a given $D$, we want primes $a,b,c,d$ with
 \be\label{eq:n2}
 ad-bc\ =\ D.
 \ee
The set of solutions in which at least one of the entries is the even prime $2$ is again of ``special form,'' and may be discarded without affecting Zariski density. Thus restricting to odd primes, we see immediately that there is a local obstruction to solving \eqref{eq:n2}, namely $D$ had better be even. 
(In fact, it is not hard to convince oneself that in the $n\times n$ case, there is again a local obstruction unless $D\equiv0(\mod 2^{n-1})$,
which is why we chose $D=4$ in \eqref{eq:soln3}.)

But now \eqref{eq:n2} looks like a ``twin prime'' type question: When can an even number $D$ be written as the difference, not of two primes, but two $E_{2}$'s? (An ``$E_{2}$'' is a number which is the product of exactly two primes.)
Miraculously, the \hyperref[GPY]{GPY} technology, extended to this setting by Goldston-Graham-Pintz-Y{\i}ld{\i}r{\i}m \cite{GoldstonGrahamPintzYildirim2009}, is able to settle the ``Bounded Gaps for $E_{2}$'s Problem,'' proving that 
$E_{2}$'s differ by at most $6$ infinitely often.
Thus there are many solutions to \eqref{eq:n2} in the primes for at least one value of $D$ in $\{2, 4,6\}$, but we do not know which!
\\

Turning now to the higher rank setting of $n\ge3$, the following clever observation was made by Nevo-Sarnak \cite{NevoSarnak2009}. One can first populate all but the last row
with primes, writing
$$
M=
\bp
*&\cdots&*&*\\
\vdots&*&*&*\\
*&*&*&*\\
m_{n,1}&\cdots&m_{n,n-1}&m_{n,n}\\
\ep
,
$$
say, where each $*$ is a prime and the $m_{n,j}$'s are variables. 
Then the equation $\det M=D$ is a {\it linear} equation to be solved in $n\ge3$ prime unknowns. For example, we found \eqref{eq:soln3} by
setting $D=4$ and
finding the solution 
$$
(a,b,c)\ =\ (5,13,19)
$$ 
to
$$
4\ =\ 
\det
\matthree
357
{11}{13}{17}
abc
\ =\ 
-6 a + 26 b - 16 c
.
$$
It goes back to I. M. Vinogradov (1937) that linear equations in at least three unknowns  can be solved in primes, and thus (overcoming many technicalities to get this simple description to actually work) Nevo-Sarnak are able to completely resolve the higher rank problem.

\subsubsection{Prime Norms in $\SL_{2}(\Z)$}\

Here is another problem of Affine Sieve type: 
Instead of restricting the entries to be prime as above, let us look at the full group $\SL_{2}(\Z)$, say, and consider its set of square-norms. 
That is, consider the set $\cS$ of values of 
$$
a^{2}+b^{2}+c^{2}+d^{2}
,
$$
where $ad-bc=1$. Does the set $\cS$ contain an infinitude of primes? 

Again one can apply the general Affine Sieve procedure, but 
Friedlander-Iwaniec 
\cite{FriedlanderIwaniec2009} found a more profitable approach. After a linear change of variables, the
 problem can 
 be converted into solving the system
\be\label{eq:FI}
\twocase{}{x^{2}+y^{2}\ =\ p+2}{}
{z^{2}+w^{2}\ =\ p-2}{}
\ee
for primes $p$ and integers $x,y,z,w$; that is, we must
write both $p+2$ and $p-2$ as sums of two squares. 
Using a ``half-dimensional'' sieve and assuming the \hyperref[conj:EH]{Elliott-Halberstam Conjecture}, 
Friedlander-Iwaniec are able to solve the system \eqref{eq:FI}, thereby (conditionally) resolving the problem in this setting.

\subsubsection{
Pseudorandom 
Primes}\label{sec:zaremba}\

The oldest (and arguably simplest)
 pseudorandom number generator is the map 
$$
x\ \mapsto\ gx(\mod p),
$$ 
where $p$ is a prime and $g$ is a primitive root mod $p$, that is, a generator of $(\Z/p\Z)^{\times}$.
For optimal equidistribution (and many other applications; see, e.g., the discussion in \cite[\S2]{Kontorovich2013}), one needs the continued fraction expansion 
\be\label{eq:cfe}
\frac gp
\ =\ 
[a_{1},a_{2},\dots,a_{k}]
\ = \ 
\cfrac{1}{a_{1}+\cfrac{1}{a_{2}+\ddots}}
\ee 
to have only ``small'' partial quotients, $a_{j}\le A$, say, for some constant $A>0$. 
Does there exist an absolute constant $A>0$ so that infinitely many such fractions $g/p$ can be found with partial quotients bounded by $A$?


To turn this into an Affine Sieve problem, observe that \eqref{eq:cfe} is equivalent to
$$
\mattwo011{a_{1}}
\mattwo011{a_{2}}
\cdots
\mattwo011{a_{k}}
\ =\
\mattwo*g*p
.
$$
Hence to find such pairs $(g,p)$, one should look at the set of second columns  in the {\it semi-group}
\be\label{eq:Gsemi}
\G\ :=\ \<\mattwo 011a:a\le A\>^{+} \cap\ \SL_{2}.
\ee
For $A\ge2$, this semigroup is Zariski dense in $\SL_{2}$, but it is known to be thin (and it is here that we wish to extend the definition of thinness beyond the realm of groups
). Instead of using the Affine Sieve, Bourgain and the author \cite{BourgainKontorovich2011, BourgainKontorovich2014} 
developed a version of
 the Hardy-Littlewood circle method to
 attack
  this problem, giving an affirmative answer to the above question: 
There are infinitely many primes $p$ and primitive roots $g(\mod p)$ so that $g/p$ has all partial quotients bounded by $A=51$.\footnote{Added in print: Shinnyih Huang \cite{Huang2013} has recently reduced this number to $A=7$, using refinements due to Frolenkov-Kan \cite{FrolenkovKan2013}.} 
In fact, they proved a ``density'' version of Zaremba's Conjecture:  Almost every natural number occurs in the set $\cS$ of bottom right entries of a matrices in $\G$ (see \cite{BourgainKontorovich2014} or \cite[\S2]{Kontorovich2013} for a precise statement). 
Thus while $\G$ is thin, the set $\cS$ is not, and no sifting is needed to produce primes in $\cS$.

\subsubsection{Prime Apollonian Curvatures}\

It seems these days no 
discourse
on
thin groups is complete without 
mention of
Apollonian gaskets. 
Lest we bore the reader, we will not yet again repeat the definitions and pictures, which are readily available elsewhere, e.g., \cite[\S3]{Kontorovich2013}. 
Nevertheless, the following question is quintessential Affine Sieve: Given a primitive Apollonian gasket $\sG$, 
which primes arise as curvatures in $\sG$?

It was proved by Sarnak \cite{SarnakToLagarias} that infinitely many prime curvatures arise, by finding primitive values of shifted binary quadratic quadratic forms among the curvatures and 
applying Iwaniec's ``half-dimensional'' sieve. 
In this way, he proved that the number of primes up to $x$ which are curvatures in $\sG$ is at least of order $x(\log x)^{-3/2}$. Bourgain \cite{Bourgain2012} sharpened the lower bound to $x(\log x)^{-1}$, that is, a positive proportion of the primes arise. Finally, Bourgain and the author \cite{BourgainKontorovich2012}, again using the circle method instead of the Affine Sieve, obtained an asymptotic formula for this number.

As in \S\ref{sec:zaremba}, this is an easy consequence of the stronger theorem that an asymptotic ``local-global'' principle holds for such curvatures (see \cite{BourgainKontorovich2012} and \cite[\S3]{Kontorovich2013} for details). So while the group and orbit in this context are again thin, the set of all curvatures is not, and the primes are obtained as  a byproduct. 

\subsection{
Improving Levels of Distribution in
the Affine Sieve}\label{sec:AffLevels}\

We conclude our discussion with
two final examples in which one can go beyond the general theory. In these, one is currently not able to produce primes, but instead can improve on the exponent of distribution over that in \eqref{eq:expDist}, without making new progress on spectral gaps as in \eqref{eq:1003ssss}. 
The idea is to avoid
putting the individual estimate \eqref{eq:rqBnd} into the sum \eqref{eq:rqErr},
and instead to try to exploit
 cancellation
  from the sum on $q$ up to $Q$, in 
  some
  analogy with the \hyperref[conj:EH]{Elliott-Halberstam Conjecture}. 
 It is not known how to do this in the general Affine Sieve, but for the specific examples 
below, 
such 
estimates
have recently been
 obtained
 by Bourgain and the author \cite{BourgainKontorovich2013, BourgainKontorovich2013a}. 
 
\subsubsection{McMullen's Arithmetic Chaos Conjecture}\

We will not describe the origins and implications of McMullen's (Classical)  Arithmetic Chaos Conjecture, referring the reader to his fascinating paper \cite{McMullen2009} and online lecture notes \cite{McMullenNotes}. The conjecture is implied by an analogue of Zaremba's Conjecture, 
purporting 
 that,
for some $A>1$, 
every sufficiently large integer arises (with the ``right'' multiplicity) in the set $\cS$ of traces of matrices in the semigroup $\G$ in  \eqref{eq:Gsemi}. At the moment, even a ``density'' version of this statement, as in \S\ref{sec:zaremba}, seems out of reach, but one can ask instead if infinitely many primes
appear in $\cS$. Not surprisingly, the standard Affine Sieve procedure applies here just as well (now requiring the work of Bourgain-Gamburd-Sarnak \cite{BourgainGamburdSarnak2011} to prove a statement functionally as strong as \eqref{eq:1003ssss}). But, if applied directly, this produces a terribly poor exponent of distribution $\vt$ in \eqref{eq:expDist}, owing to the terribly poor ``spectral gap'' $\gT$. 
Using different tools in this setting
, Bourgain and the author \cite{BourgainKontorovich2013a} have produced in this context an unconditional exponent of distribution $\vt=1/4$, thus showing that $\cS$ contains $R$-almost-primes with $R=5$.

\subsubsection{Thin Pythagorean Hypotenuses}\

Finally, let us return 
again to the Pythagorean setting 
of the quadratic form $F$ in \eqref{eq:FPyth}, the cone $F=0$, the base point $\bx_{0}=(3,4,5)$, and a thin group $\G$ as in \S\ref{sec:thin}. Instead of studying areas, let us now take as our function $f$ the ``hypotenuse,'' $f(\bx)=z$. 
Do infinitely many primes arise in $\cS=f(\cO)$?

If $\cO$ were the full orbit of all Pythagorean triples, then, through the parametrization \eqref{eq:PythParam}, we would essentially asking whether primes can be represented as sums of two squares. 
As is very well-known,
 Fermat answered in the affirmative 
almost 
400 years ago, namely all primes $\equiv1(\mod4)$ are hypotenuses.

But in the thin setting, it seems quite difficult to
produce primes at this time.
%
One new difficulty 
here
is that, unlike 
other
problems
 described above,
we now have
 not only a thin orbit $\cO$, but the set $\cS$ of hypotenuses is itself thin! 
The number of integers in $\cS$ up to $x$, even with multiplicity, is about $x^{\gd}$, where $\gd<1$ is the growth exponent of $\G$ as in \eqref{eq:gdIs}. 

What does the Affine Sieve process give? Returning to the exponent of distribution $\vt$ in \eqref{eq:expDist}, we see that the degree of the hypotenuse function is $d=1$, and the value of $C$ can 
be whittled down to $2$ as  in \S\ref{sec:PythA}. Moreover, to try to optimize $\vt$, we can restrict our attention to thin groups $\G$ whose growth exponent $\gd$ is almost as large as possible, $\gd=1-\vep$.
Then, even assuming a ``square-root'' version of \eqref{eq:1003ssss}, that is, assuming the ``spectral gap'' can be set to $\gT=1/2$, 
we obtain a (very conditional) exponent of distribution $\vt= 1/4-\vep$, producing $R$-almost-primes in $\cS$ with $R=5$. 
In
\cite{BourgainKontorovich2013},
Bourgain and the author
obtained, again for $\G$ having growth exponent $\gd$ sufficiently close to $1$, the exponent of distribution $\vt=7/24-\vep$  unconditionally, thereby producing $R$-almost-primes with $R=4$ in this thin setting. 
The methods (bilinear forms, exponential sums, and dispersion) are outside the scope of this survey. 
\\

\

{\bf Added in proof:}
The explicit values of $R$ in \S\S\ref{sec:genApply}--\ref{sec:Thin} are  now outdated; recent work of the author and Jiuzu Hong \cite{HongKontorovich2014} gives an improvement on the general Affine Sieve procedure which differs slightly from that given here (we will not go into the technicalities). Still, the problem of going beyond these values in specific cases remains, and in these settings (that is, in \S\ref{sec:AffLevels}), the reported $R$ values are still the best known.

\

\




\section*{Acknowledgements}\

The author is grateful to Dick Gross for the invitation to visit Harvard University, during which time  these lecture notes were written. 
Many thanks to
 John Friedlander,
 Andrew Granville,
 Curt McMullen,
 Sam Payne,
 Peter Sarnak,
 Yitang Zhang, and the referee for comments and suggestions on an earlier draft.
Thanks also to the organizers of the ``Hyperbolic Geometry and Arithmetic'' 
workshop
in Toulouse,  November 2012, especially Cyril Lecuire;  the second half of these notes 
was
conceived at this 
meeting
(the first half was not yet a theorem).
The author is also indebted  to Jean-Pierre Otal, without whose persistence these notes would not have materialized.

\newpage

\bibliographystyle{alpha}

\bibliography{../../../AKbibliog}

\end{document}